\documentclass[12pt,fleqn, reqno]{article}

\usepackage{mathrsfs}
\usepackage{amsfonts}
\usepackage{amsthm}
\usepackage{amsmath}
\usepackage{graphicx}
\usepackage{amscd,amssymb,amsthm, color}
\usepackage{setspace}
\numberwithin{equation}{section}

\begin{document}

\allowdisplaybreaks

\thispagestyle{plain}
\vspace{4cc}
\begin{center}
{\large \bf Further Novel reductions of Kamp\'e de F\'eriet function}\\\vspace{2cc}
{\bf Arjun Kumar Rathie}\\
\medskip
 Department of Mathematics, \\
 Vedant College of Engineering and Technology\\ (Rajasthan Technical University), Bundi-323021,\\ Rajasthan, India\\
 E-mail: {\tt arjunkumarrathie@gmail.com}
 \end{center}
\vspace{2cc}
 {\small {\bf Abstract}. In a recent paper, Rathie and Pogany established thirty two novel and general reductions of two and three variables generalized hypergeometric functions. In this paper we provide twenty four further novel and general reduction formulas. The results are established by the application of Beta and Gamma integral methods to the three identities involving  
products of 
generalized hypergeometric functions obtained earlier by Kim and Rathie. As special cases, we mention some interesting results. 
}\\

\noindent {\bf 2010 Mathematics Subject Classification :}\\ Primary: 33C65,~ 33C20; ~\\
Secondary:  33B15,~ 33C05,~ 33C15 \\

\noindent {\bf Keywords and Phrases}. Generalized hypergeometric function, Double and triple  hypergeometric functions. 
\section{Main Reduction Formulas} 
The results to be proved in this paper are given in the following three theorems.\\
 
 \noindent
\textbf{Theorem 1 : } 
For $m, n \in \mathbb N_0$,   $\Re(e)>\Re(d)>0$ and $\alpha, \beta \in \mathbb C$, the following results hold true.   \small \begin{align}  
 & F_{1:\,1;1}^{1:\,1;1}    \left[\begin{array}{c} d:\,\alpha\,;\, \beta   \\  e : 2\alpha+m; 2\beta+n \end{array} \left| \begin{array}{c} x\\x \end{array}  \right.\right] \nonumber \\
 ~~& = \sum_{j  = 0}^m \sum_{k = 0}^n \frac{(d)_{j+k}(-m)_j (-n)_k (2\alpha-1)_j (2\beta-1)_k 	   x^{j+k}}{(e)_{j+k}(2\alpha +m)_j  (2\beta+n)_k \left(\alpha-\frac{1}{2}\right)_j \left(\beta-\frac{1}{2}\right)_k\, 2^{2j+2k}~ j!~ k!} \nonumber \\
&\times \mathcal{S}_{1: 0 ; 3}^{1:\,0; 2}\left[\begin{array}{c}
{[d+j+k: 1,2]:-; \frac{ (\alpha+\beta+j+k)}{2},~  \frac{ (\alpha+\beta+j+k+1)}{2}} \\ {[e+j+k: 1,2]:-; \alpha+j+\frac{1}{2},~  \beta+k+\frac{1}{2},~ \alpha+\beta+j+k}\end{array} \left| \begin{array}{c} x \\ \frac{1}{4} x^{2} \end{array}\right. \right]  
\end{align}   
\small
 \begin{align} 
 & F_{1:~1;1}^{ 1:~1;1}  \left[\begin{array}{c} d:\,  \alpha; \, \beta +n \\    e : 2\alpha+m; 2\beta + n \end{array}\left|\begin{array}{c} x\\-x \end{array} \right. \right] \nonumber \\
~~& = F^{(3)}\left[\begin{array}{c} d::\, -\,;\, -\,;\, -\,;\, -\,;\, \alpha\,;\,\beta\,\\  e:: -;-;-;-; 2\alpha+m; 2\beta+n \end{array} \left| -x, x, x \right. \right] \notag \\
&  = \sum_{j  = 0}^m \sum_{k = 0}^n  \frac{(d)_{j+k}~(-m)_j~ (-n)_k (2\alpha-1)_j (2\beta-1)_k\, x^{j+k}}{(e)_{j+k} (2\alpha+m)_j
(2\beta+n)_k  \left(\alpha-\frac{1}{2}\right)_j~\left(\beta-\frac{1}{2}\right)_k ~ 2^{2j+2k}\,j!\,k!} \notag \\
 &  ~~\times   {}_4F_5 \left[ \begin{array}{c} \frac{  \alpha + \beta + j+k }{2},~ \frac{ \alpha + \beta + j+k+1 }{2}, \frac{d+j+k}{2}, \, \frac{ d+j+k+1}{2} \\ \alpha + j + \frac{1}{2}, \beta + k + \frac{1}{2},\, \alpha+\beta+j+k,\,\frac{ e+j+k}{2},\, \frac{ e+j+k+1}{2}\end{array} \left| \frac{1}{4}{x^2} \right. \right]
\end{align}
\small \begin{align} 
& F_{0:\,1;1}^{1:\,1;1}  \left[\begin{array}{c} d:\,\alpha\,;\, \beta   \\ - : 2\alpha+m; 2\beta+n \end{array} \left| \begin{array}{c} x\\x \end{array}  \right.\right] \nonumber \\
 ~~& = \sum_{j  = 0}^m \sum_{k = 0}^n \frac{(d)_{j+k}(-m)_j (-n)_k (2\alpha-1)_j (2\beta-1)_k 	   x^{j+k}}{ (2\alpha +m)_j  (2\beta+n)_k \left(\alpha-\frac{1}{2}\right)_j \left(\beta-\frac{1}{2}\right)_k\, 2^{2j+2k}~ j!~ k!} \nonumber \\
&\times \mathcal{S}_{0: 0 ; 3}^{1:\,0; 2}\left[\begin{array}{c}
{[d+j+k: 1,2]:-; \frac{ \alpha+\beta+j+k}{2}  ; \frac{  \alpha+\beta+j+k+1}{2}} \\{-:~ -;~   \alpha+j+\frac{1}{2}; \beta+k+\frac{1}{2};  \alpha+\beta+j+k}\end{array} \left| \begin{array}{c} x \\ \frac{1}{4} x^{2} \end{array}\right. \right]  
\end{align}   
\small \begin{align} 
 &F_{0:~1;1}^{ 1:~1;1} \left[\begin{array}{c} d:\,  \alpha; \, \beta +n \\    - : 2\alpha+m; 2\beta + n \end{array}\left|\begin{array}{c} x\\-x \end{array} \right. \right] \nonumber \\
~~& = F^{(3)}\left[\begin{array}{c} d::\, -\,;\, -\,;\, -\,;\, -\,;\, \alpha\,;\,\beta\,\\ -:: -;-;-;-; 2\alpha+m; 2\beta+n \end{array} \left| -x, x, x \right. \right] \notag \\
&  = \sum_{j  = 0}^m \sum_{k = 0}^n  \frac{(d)_{j+k}~(-m)_j~ (-n)_k (2\alpha-1)_j (2\beta-1)_k\, x^{j+k}}{  (2\alpha+m)_j
(2\beta+n)_k  \left(\alpha-\frac{1}{2}\right)_j~\left(\beta-\frac{1}{2}\right)_k ~ 2^{2j+2k}\,j!\,k!} \notag \\
 &  \times   {}_4F_3 \left[ \begin{array}{c} \frac{  \alpha + \beta + j+k}{2},~ \frac{ \alpha + \beta + j+k+1 }{2}, \frac{  d+j+k}{2}, \, \frac{  d+j+k+1}{2} \\ \alpha + j + \frac{1}{2}, \beta + k + \frac{1}{2},\, \alpha+\beta+j+k\end{array} \left| \frac{1}{4}{x^2} \right. \right]
\end{align}

\noindent
\textbf{Theorem 2 : } 
For $m, n \in \mathbb N_0$,   $\Re(e)>\Re(d)>0$ and $\alpha, \beta \in \mathbb C$, the following results hold true.   \small \begin{align}  
& F_{1:\,1;1}^{1:\,1;1}  \left[\begin{array}{c} d:\,\alpha\,;\, \beta   \\  e : 2\alpha-m; 2\beta-n \end{array} \left| \begin{array}{c} x\\x \end{array}  \right.\right] \nonumber \\
 ~~& = \sum_{j  = 0}^m \sum_{k = 0}^n \frac{(d)_{j+k}\, (-1)^{j+k}~(-m)_j (-n)_k (2\alpha-2m-1)_j (2\beta-2n-1)_k   x^{j+k}}{(e)_{j+k}(2\alpha -m)_j  (2\beta-n)_k \left(\alpha-m-\frac{1}{2}\right)_j \left(\beta-n-\frac{1}{2}\right)_k\, 2^{2j+2k}~ j!~ k!} \nonumber \\
&\times \mathcal{S}_{1: 0 ; 3}^{1:\,0; 2}\left[\begin{array}{c}
{[d+j+k: 1,2]:-; \frac{ \alpha+\beta+j+k-m-n}{2}  ; \frac{\alpha+\beta+j+k-m-n+1}{2}} \\ {[e+j+k: 1,2]:-; \alpha+j-m+\frac{1}{2}; \beta+k-n+\frac{1}{2};  \alpha+\beta+j+k-m-n}\end{array} \left| \begin{array}{c} x \\ \frac{1}{4} x^{2} \end{array}\right. \right]  
\end{align}   
\small \begin{align} 
 &F_{1:~1;1}^{ 1:~1;1} \left[\begin{array}{c} d:\,  \alpha; \, \beta -n \\    e : 2\alpha-m; 2\beta - n \end{array}\left|\begin{array}{c} x\\-x \end{array} \right. \right] \nonumber \\
~~& = F^{(3)}\left[\begin{array}{c} d::\, -\,;\, -\,;\, -\,;\, -\,;\, \alpha\,;\,\beta\,\\  e:: -;-;-;-; 2\alpha-m; 2\beta-n \end{array} \left| -x, x, x \right. \right] \notag \\
&  = \sum_{j  = 0}^m \sum_{k = 0}^n  \frac{(d)_{j+k}~(-1)^{j+k}~(-m)_j~ (-n)_k (2\alpha-2m-1)_j (2\beta-2n-1)_k\, x^{j+k}}{(e)_{j+k} (2\alpha-m)_j
(2\beta-n)_k  \left(\alpha-m-\frac{1}{2}\right)_j~\left(\beta-n-\frac{1}{2}\right)_k ~ 2^{2j+2k}\,j!\,k!} \notag \\
 &  \times   {}_4F_5 \left[ \begin{array}{c} \frac{  \alpha + \beta + j+k-m-n }{2},~ \frac{ \alpha + \beta + j+k-m-n+1 }{2}, \frac{ d+j+k}{2}, \, \frac{ d+j+k+1}{2} \\ \alpha + j-m + \frac{1}{2}, \beta + k-n + \frac{1}{2},\, \alpha+\beta+j+k-m-n,\,\frac{ e+j+k}{2},\, \frac{ e+j+k+1}{2}\end{array} \left| \frac{1}{4}{x^2} \right. \right]
\end{align}
\small \begin{align} 
& F_{0:\,1;1}^{1:\,1;1}  \left[\begin{array}{c} d:\,\alpha\,;\, \beta   \\ - : 2\alpha-m; 2\beta-n \end{array} \left| \begin{array}{c} x\\x \end{array}  \right.\right] \nonumber \\
 ~~& = \sum_{j  = 0}^m \sum_{k = 0}^n \frac{(d)_{j+k}\, (-1)^{j+k}~(-m)_j (-n)_k (2\alpha-2m-1)_j (2\beta-2n-1)_k   x^{j+k}}{ (2\alpha -m)_j  (2\beta-n)_k \left(\alpha-m-\frac{1}{2}\right)_j \left(\beta-n-\frac{1}{2}\right)_k\, 2^{2j+2k}~ j!~ k!} \nonumber \\
&\times \mathcal{S}_{0: 0 ; 3}^{1:\,0; 2}\left[\begin{array}{c}
{[d+j+k: 1,2]:-; \frac{ \alpha+\beta+j+k-m-n}{2}  ; \frac{ \alpha+\beta+j+k-m-n+1}{2}} \\ {-:-; \alpha+j-m+\frac{1}{2}; \beta+k-n+\frac{1}{2};  \alpha+\beta+j+k-m-n}\end{array} \left| \begin{array}{c} x \\ \frac{1}{4} x^{2} \end{array}\right. \right]  
\end{align}   
\small \begin{align} 
 &F_{0:~1;1}^{ 1:~1;1} \left[\begin{array}{c} d:\,  \alpha; \, \beta-n  \\  - : 2\alpha-m; 2\beta - n \end{array}\left|\begin{array}{c} x\\-x \end{array} \right. \right] \nonumber \\
~~& = F^{(3)}\left[\begin{array}{c} d::\, -\,;\, -\,;\, -\,;\, -\,;\, \alpha\,;\,\beta\,\\  -:: -;-;-;-; 2\alpha-m; 2\beta-n \end{array} \left| -x, x, x \right. \right] \notag \\
&  = \sum_{j  = 0}^m \sum_{k = 0}^n  \frac{(d)_{j+k}~(-1)^{j+k}~(-m)_j~ (-n)_k (2\alpha-2m-1)_j (2\beta-2n-1)_k\, x^{j+k}}{ (2\alpha-m)_j
(2\beta-n)_k  \left(\alpha-m-\frac{1}{2}\right)_j~\left(\beta-n-\frac{1}{2}\right)_k ~ 2^{2j+2k}\,j!\,k!} \notag \\
 &  \times   {}_4F_3 \left[ \begin{array}{c} \frac{  \alpha + \beta + j+k-m-n }{2},~ \frac{ \alpha + \beta + j+k-m-n+1 }{2}, \frac{ d+j+k}{2}, \, \frac{ d+j+k+1}{2} \\ \alpha + j-m + \frac{1}{2}, \beta + k-n + \frac{1}{2},\, \alpha+\beta+j+k-m-n \end{array} \left| \frac{1}{4}{x^2} \right. \right]
\end{align}
\noindent \textbf{Theorem 3 :} For $m, n \in \mathbb N_0$,   $\Re(e)>\Re(d)>0$ and $\alpha, \beta \in \mathbb C$, the following results hold true.  
\small \begin{align}  
& F_{1:\,1;1}^{1:\,1;1}  \left[\begin{array}{c} d:\,\alpha\,;\, \beta   \\  e : 2\alpha+m; 2\beta-n \end{array} \left| \begin{array}{c} x\\x \end{array}  \right.\right] \nonumber \\
 ~~& = \sum_{j  = 0}^m \sum_{k = 0}^n \frac{(d)_{j+k}~(-1)^k ~(-m)_j (-n)_k (2\alpha-1)_j (2\beta-2n-1)_k 	   x^{j+k}}{(e)_{j+k}(2\alpha +m)_j  (2\beta-n)_k \left(\alpha-\frac{1}{2}\right)_j \left(\beta-n-\frac{1}{2}\right)_k\, 2^{2j+2k}~ j!~ k!} \nonumber \\
&\times \mathcal{S}_{1: 0 ; 3}^{1:\,0; 2}\left[\begin{array}{c}
{[d+j+k: 1,2]:-; \frac{ \alpha+\beta+j+k-n}{2}  ; \frac{ \alpha+\beta+j+k-n+1}{2}} \\ {[e+j+k: 1,2]:-; \alpha+j+\frac{1}{2}; \beta+k-n+\frac{1}{2};  \alpha+\beta+j+k-n}\end{array} \left| \begin{array}{c} x \\ \frac{1}{4} x^{2} \end{array}\right. \right]  
\end{align}   
\small \begin{align} 
 &F_{1:~1;1}^{ 1:~1;1} \left[\begin{array}{c} d:\,  \alpha; \, \beta-n \\    e : 2\alpha+m; 2\beta - n \end{array}\left|\begin{array}{c} x\\-x \end{array} \right. \right] \nonumber \\
~~& = F^{(3)}\left[\begin{array}{c} d::\, -\,;\, -\,;\, -\,;\, -\,;\, \alpha\,;\,\beta\,\\  e:: -;-;-;-; 2\alpha+m; 2\beta-n \end{array} \left| -x, x, x \right. \right] \notag \\
&  = \sum_{j  = 0}^m \sum_{k = 0}^n  \frac{(d)_{j+k}~(-1)^k~(-m)_j~ (-n)_k (2\alpha-1)_j (2\beta-2n-1)_k\, x^{j+k}}{(e)_{j+k} (2\alpha+m)_j
(2\beta-n)_k  \left(\alpha-\frac{1}{2}\right)_j~\left(\beta-n-\frac{1}{2}\right)_k ~ 2^{2j+2k}\,j!\,k!} \notag \\
 &  \times   {}_4F_5 \left[ \begin{array}{c} \frac{  \alpha + \beta + j+k-n}{2},~ \frac{ \alpha + \beta + j+k-n+1 }{2}, \frac{1 d+j+k}{2}, \, \frac{ d+j+k+1}{2} \\ \alpha + j + \frac{1}{2}, \beta + k -n+ \frac{1}{2},\, \alpha+\beta+j+k-n,\,\frac{ e+j+k}{2},\, \frac{ e+j+k+1}{2}\end{array} \left| \frac{1}{4}{x^2} \right. \right]
\end{align}
\small \begin{align} 
& F_{0:\,1;1}^{1:\,1;1}  \left[\begin{array}{c} d:\,\alpha\,;\, \beta   \\ - : 2\alpha+m; 2\beta-n \end{array} \left| \begin{array}{c} x\\x \end{array}  \right.\right] \nonumber \\
 ~~& = \sum_{j  = 0}^m \sum_{k = 0}^n \frac{(d)_{j+k}~(-1)^k~(-m)_j (-n)_k (2\alpha-1)_j (2\beta-2n-1)_k 	   x^{j+k}}{(e)_{j+k}~ (2\alpha +m)_j  (2\beta-n)_k \left(\alpha-\frac{1}{2}\right)_j \left(\beta-n-\frac{1}{2}\right)_k\, 2^{2j+2k}~ j!~ k!} \nonumber \\
&\times \mathcal{S}_{0: 0 ; 3}^{1:\,0; 2}\left[\begin{array}{c}
{[d+j+k: 1,2]:-; \frac{ \alpha+\beta+j+k-n}{2}  ; \frac{ \alpha+\beta+j+k-n+1}{2}} \\{-:~ -;~   \alpha+j+\frac{1}{2}; \beta+k-n+\frac{1}{2};  \alpha+\beta+j+k-n}\end{array} \left| \begin{array}{c} x \\ \frac{1}{4} x^{2} \end{array}\right. \right]  
\end{align}   
\small \begin{align} 
 &F_{0:~1;1}^{ 1:~1;1} \left[\begin{array}{c} d:\,  \alpha; \, \beta -n \\    - : 2\alpha+m; 2\beta - n \end{array}\left|\begin{array}{c} x\\-x \end{array} \right. \right] \nonumber \\
~~& = F^{(3)}\left[\begin{array}{c} d::\, -\,;\, -\,;\, -\,;\, -\,;\, \alpha\,;\,\beta\,\\ -:: -;-;-;-; 2\alpha+m; 2\beta-n \end{array} \left| -x, x, x \right. \right] \notag \\
&  = \sum_{j  = 0}^m \sum_{k = 0}^n  \frac{(d)_{j+k}~(-1)^k~(-m)_j~ (-n)_k (2\alpha-1)_j (2\beta-2n-1)_k\, x^{j+k}}{  (2\alpha+m)_j
(2\beta-n)_k  \left(\alpha-\frac{1}{2}\right)_j~\left(\beta-n-\frac{1}{2}\right)_k ~ 2^{2j+2k}\,j!\,k!} \notag \\
 &  \times   {}_4F_3 \left[ \begin{array}{c} \frac{  \alpha + \beta + j+k-n }{2},~ \frac{ \alpha + \beta + j+k-n+1 }{2}, \frac{ d+j+k}{2}, \, \frac{ d+j+k+1}{2} \\ \alpha + j + \frac{1}{2}, \beta + k -n+ \frac{1}{2},\, \alpha+\beta+j+k-n\end{array} \left| \frac{1}{4}{x^2} \right. \right]
\end{align}
Proofs : The proofs of the results given in the Theorems 1 , 2 and 3 are same as given in the paper[3]. However, while deriving the results, we shall make use of the results (2.1), (2.2) and (2.3) given in the paper [1].\\

The details are given in [2].
\section{Special Cases}
\noindent 1. In theorems 1, 2 and 3, if we set $\beta=\alpha$, we get the known results recently obtained in [3].\\
2. In theorems 1 or 2 or 3, if we set $m=n=0$, we get the following interesting results.
\begin{align}
 F_{1:\,1;1}^{1:\,1;1}&  \left[\begin{array}{c} d:\,\alpha\,;\, \beta   \\  e : 2\alpha ; 2\beta  \end{array} \left| \begin{array}{c} x\\x \end{array}  \right.\right] \nonumber \\
  ~~ &= \mathcal{S}_{1: 0 ; 3}^{1:\,0; 2}\left[\begin{array}{c}
{[d : 1,2]:-; \frac{1}{2}(\alpha+\beta ), \frac{1}{2}(\alpha+\beta +1)} \\ {[e : 1,2]:-; \alpha +\frac{1}{2}, \beta+ \frac{1}{2},  \alpha+\beta }\end{array} \left| \begin{array}{c} x \\ \frac{1}{4} x^{2} \end{array}\right. \right]  
\end{align} 
\begin{align} 
  F_{1:~1;1}^{ 1:~1;1}& \left[\begin{array}{c} d:\,  \alpha; \, \beta  \\    e : 2\alpha ; 2\beta   \end{array}\left|\begin{array}{c} x\\-x \end{array} \right. \right] \nonumber \\
~~& = F^{(3)}\left[\begin{array}{c} d::\, -\,;\, -\,;\, -\,;\, -\,;\, \alpha\,;\,\beta\,\\  e:: -;-;-;-; 2\alpha ; 2\beta  \end{array} \left| -x, x, x \right. \right] \notag \\
&  =    {}_4F_5 \left[ \begin{array}{c} \frac{1}{2}( \alpha + \beta   ),~ \frac{1}{2}(\alpha + \beta  +1 ), \frac{1}{2}(d ), \, \frac{1}{2}(d +1) \\ \alpha   + \frac{1}{2}, \beta   + \frac{1}{2},\, \alpha+\beta,\,\frac{1}{2}(e ),\, \frac{1}{2}(e +1)\end{array} \left| \frac{1}{4}{x^2} \right. \right]
\end{align}
\begin{align}  
  F_{0:\,1;1}^{1:\,1;1}&  \left[\begin{array}{c} d:\,\alpha\,;\, \beta   \\ - : 2\alpha ; 2\beta  \end{array} \left| \begin{array}{c} x\\x \end{array}  \right.\right] \nonumber \\
 ~~& =  \mathcal{S}_{0: 0 ; 3}^{1:\,0; 2}\left[\begin{array}{c}
{[d : 1,2]:-; \frac{1}{2}(\alpha+\beta )  ; \frac{1}{2}(\alpha+\beta +1)} \\{-:~ -;~   \alpha +\frac{1}{2}; \beta +\frac{1}{2};  \alpha+\beta }\end{array} \left| \begin{array}{c} x \\ \frac{1}{4} x^{2} \end{array}\right. \right]  
\end{align}   
\begin{align} 
 F_{0:~1;1}^{ 1:~1;1}& \left[\begin{array}{c} d:\,  \alpha; \, \beta   \\    - : 2\alpha ; 2\beta   \end{array}\left|\begin{array}{c} x\\-x \end{array} \right. \right] \nonumber \\
~~& = F^{(3)}\left[\begin{array}{c} d::\, -\,;\, -\,;\, -\,;\, -\,;\, \alpha\,;\,\beta\,\\ -:: -;-;-;-; 2\alpha ; 2\beta  \end{array} \left| -x, x, x \right. \right] \notag \\
&  =     {}_4F_3 \left[ \begin{array}{c} \frac{1}{2}( \alpha + \beta   ),~ \frac{1}{2}(\alpha + \beta  +1 ), \frac{1}{2}(d ), \, \frac{1}{2}(d +1) \\ \alpha   + \frac{1}{2}, \beta   + \frac{1}{2},\, \alpha+\beta \end{array} \left| \frac{1}{4}{x^2} \right. \right]
\end{align}
\noindent Similarly other results can be obtained.

	 \end{document}